\documentclass[12pt,a4paper]{article}
\usepackage{latexsym}
\usepackage{amsfonts}
\usepackage{amssymb}
\usepackage{amsbsy}
\usepackage{mathrsfs}
\usepackage{amsmath}
\usepackage{epsfig}
\usepackage{theorem}
\usepackage{color}
\usepackage[all]{xy}
\xyoption{all}

\def\scM{{\mathcal{M}}}
\def\Top{\mathcal{T}\!op}
\def\id{\mathop{\rm id}}
\def\proj{\mathop{\rm proj}}

\newenvironment{proof}{\par\noindent\textbf{Beweis:}}{\hfill\ensuremath{\square}}
\newenvironment{defi}{\par\noindent\textbf{Definition:}}{}
\newenvironment{theo}{\par\noindent\textbf{Theorem:}}{}
\newenvironment{prop}{\par\noindent\textbf{Proposition:}}{}

\begin{document}
\title{The HELP-Lemma And Its Converse In Quillen Model Categories}
\author{R.M.~Vogt}
\maketitle

\paragraph*{Abstract}
We show that a map $p:X\to Y$ between fibrant objects in a closed model category is a weak equivalence 
if and only
if it has the right homotopy extension lifting property with respect to all cofibrations.
The dual statement holds for maps between cofibrant objects. 

\vspace{2ex}
The HELP-Lemma states that a homotopy equivalence $p:X\to Y$ of topological spaces has the homotopy
extension lifting property , HELP for short, for all closed cofibrations \cite[Appendix Thm. 3.5]{BV}. 
The lemma, variants of it and their Eckmann-Hilton duals (e.g. see \cite[II.1.11]{Bau}, 
\cite[Thm. 4, Thm. 4*]{May}) have proven to be 
very useful tools in
homotopy theory. Surprisingly it has been overlooked that the homotopy extension lifting property
 for all closed
cofibrations characterises homotopy equivalences:

\begin{prop}\label{1} A map of topological spaces $X\to Y$ is a homotopy equivalence if
and only if it has the HELP for all closed cofibrations.
\end{prop}

In this note we will prove this statement and its Eckmann-Hilton dual in arbitrary closed model categories 
$\scM$ in the sense of Quillen \cite{Qu}, (see also \cite{DS}). Applying our result to
the category $\Top$ of toplogical spaces with the Str{\o}m model structure \cite{Str} we retrieve 
the Proposition.

\begin{defi} Let $i:A\to B$ and $p:X\to Y$ be maps in $\scM$.

(1) We say that $p$ has the right HELP with respect to $i$, if for each
 not necessarily commutative square
$$
\xymatrix{
A \ar[r]^{f_A}\ar[d]_i & X\ar[d]^p
\\
B  \ar[r]^g &Y 
}
\eqno(\ast)
$$
and each right homotopy $h_A: A\to Y^I$ from $p\circ f_A$ to $g\circ i$, where
$Y^I$ is a path object $\xymatrix{Y\ar[r]^j_\sim & Y^I\ar@{->>}[r]^\pi &Y\times Y}$ for $Y$,
there is a map $f:B\to X$ and a right homotopy $h:B\to Y^I$ from $p\circ f$ to
$g$ such that $f\circ i=f_A$ and $h\circ i=h_A$.

(2) We say that $i$ has the left HELP with respect to $p$, if for each
 not necessarily commutative square
$$
\xymatrix{
A \ar[r]^{f}\ar[d]_i & X\ar[d]^p
\\
B  \ar[r]^{g_Y} &Y 
}
$$
and each left homotopy $h_Y: Z_A\to Y$ from $g_Y\circ i$ to $p\circ f$,
where $Z_A$ is a cylinder object  $\xymatrix{A\sqcup A\;\ar@{{>}->}[r]^j & Z_A\ar[r]^\sigma_\sim & A}$
for $A$, there is a map $g:B\to X$ and a left homotopy $h:Z_A\to X$ from $g\circ i$
to $f$ such that $p\circ g=g_Y$ and $p\circ h=h_Y$.
\end{defi}

\begin{theo} (1) A map $p:X\to Y$ of fibrant objects is a weak equivalence in $\scM$
if and only if it has the right HELP with respect to all cofibrations.

(2) A map $i: A\to B$ of cofibrant objects is a weak equivalence in $\scM$
if and only if it has the left HELP with respect to all fibrations.
\end{theo}
\begin{proof}
The two statements are dual so we just prove the first one. 

Since $X$ and $Y$ are fibrant the projections
$$
\xymatrix{
X & X\times Y \ar[l]_{p_X}\ar[r]^{p_Y} & Y &\quad \textrm{and}\quad &
Y & Y\times Y \ar[l]_{p_1}\ar[r]^{p_2} & Y 
}
$$
are fibrations.

Suppose that $p$ is a weak equivalence and that we are given a square $(\ast)$. Consider the 
commutative diagram
$$
\xymatrix@R=12ex@C=8ex{
X\ar[rd]^r\ar[rddd]_{\id} \ar[rrr]^p
&&& Y\ar[rdd]^{\id}\ar[d]^\sim_j
\\
& P \ar[rr]^{q_2}\ar@{->>}[d]^{q_1} \ar @{} [drr] | {pullback}
&& Y^I\ar@{->>}[d]^{\pi}\\
& X\times Y \ar[rr]^{p\times\id}\ar@{->>}[d]^{\proj_X} \ar @{} [drr] | {pullback}
&& Y\times Y\ar@{->>}[r]^{\proj_2} \ar@{->>}[d]^{\proj_1}
& Y
\\
& X \ar[rr]^{p}
&& Y
}
$$
where $r=(\id,j\circ p)$. Since $\proj_1\circ\pi$ is a weak equivalence, so is $\proj_X\circ q_1$ and hence
 $r$. It follows
that $q_2$ and $\proj_2\circ\pi\circ q_2$ 
are weak equivalences. Hence $q_3=\proj_Y\circ q_1=\proj_2
\circ\pi\circ q_2:P\twoheadrightarrow Y$ is a weak equivalence. Now consider
$$
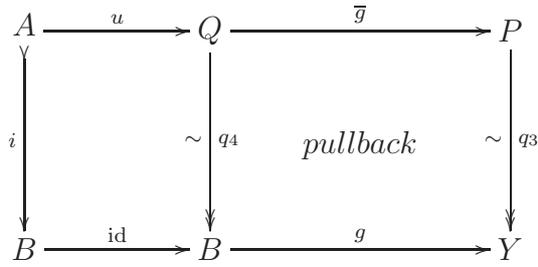
\xymatrix{
\stackrel{\displaystyle{A}}{\phantom{,}}\ar[rr]^u \ar@{>->}[dd]_i && Q\ar[rr]^{\overline{g}}
 \ar@{->>}[dd]^{q_4}_\sim
 && P\ar@{->>}[dd]^{q_3}_\sim
\\
&&&{pullback}&\\
B \ar[rr]^{\id} && B\ar[rr]^g && Y
}
$$
where $u=(i,v)$ and $v:A\to P$ is induced by $(f_A,g\circ i):A\to X\times Y$ and $h_A:A\to Y^I$. Since $q_4$ is
a trivial fibration there is a section $s:B\to Q$ such that $s\circ i=u$. Define
$$
\begin{array}{l}
f =\proj_X\circ q_1\circ\overline{g}\circ s: B\to X\\
h: q_2\circ\overline{g}\circ s: B\to Y^I
\end{array}
$$

Conversely, suppose that $p$ has the right HELP.  
Consider the diagram
$$
\xymatrix{
\stackrel{\displaystyle{A}}{\phantom{,}}\ar@{>->}[dd]^i \ar[rr]^u
&& P \ar[rr]^{\proj_X\circ q_1} \ar@{->>}[rrdd]_{q_3}
&& X\ar[dd]^p
\\
& I &&&\\
B \ar[rrrr]_g  &&&& Y
}
$$
where square $I$ is supposed to commute. We define
$$
f_A={\proj}_X\circ q_1\circ u:A\to X \quad\textrm{ and }\quad h_A=q_2\circ u:A\to Y^I.
$$
Then $h_A$ is a right homotopy from $p\circ f_A$ to $g\circ i$. Hence there exist
$$
f:B\to X \quad \textrm{ and }\quad h:B\to Y^I
$$
such that $h$ is a right homotopy from $p\circ f$ to $g$ and $f\circ i=f_A$ and $h\circ i=h_A$.
Then $f$ and $h$ induce a map
$$
k:B\to P
$$
such that $k\circ i=u$ and $q_3\circ k=\proj_2\circ\pi\circ k=g$.

Hence $q_3$ has the right lifting property with respect to all cofibration and has to be a trivial fibration.
Since $q_2:P\to Y^I$ is a right homotopy from $p\circ\proj_X\circ q_1$ to $q_3$ and since a
 map right homotopic 
to a weak equivalence is itself a weak equivalence, $p\circ\proj_X\circ q_1$ is a weak equivalence.
 Since $\proj_X\circ
q_1$ is a weak equivalence, so it $p$.
\end{proof}


\begin{thebibliography}{99}
\bibitem{Bau}  H.J. Baues, Algebraic homotopy, Cambridge University Press 1989.

\bibitem{BV} J.M. Boardman, R.M. Vogt, 
Homotopy invariant structures on topological spaces,
Lecture Notes in Math. 347, Springer Verlag, Berlin 1973.

\bibitem{DS}  
W.G. Dwyer, J. Spalinski, Homotopy theories and model categories, Handbook of Algebraic Topology
(I.M. James, ed.), Elsevier Science B.V., 1995.

\bibitem{May} J.P. May, The dual Whitehead theorems, Topological Topics (I.M. James, ed.),London
Math. Soc. Lecture Notes Ser. 86, Cambridge University Press 1983. 

\bibitem{Qu} D.G. Quillen, Rational homotopy theory, Ann. Math. 90 (1969), 205-295.

\bibitem{Str} A. Str{\o}m, The homotopy category is a homotopy 
category, Arch. Math. 23 (1972), 435-441.




\end{thebibliography}
\end{document}